\documentclass[pdflatex,sn-basic,Numbered]{sn-jnl}
\usepackage[T1]{fontenc}
\usepackage{graphicx}
\DeclareGraphicsExtensions{.pdf,.png,.jpg}

\usepackage{amsmath,amssymb,amsfonts}
\usepackage{tikz}
\usepackage{float}
\usepackage{hyperref}
\usepackage{cleveref}
\usepackage{tikz-network}
\usetikzlibrary{positioning,quotes}
\usepackage{comment}

\newtheorem{theorem}{Theorem}
\newtheorem{lemma}{Lemma}
\newtheorem{claim}{Claim}
\newtheorem{corollary}{Corollary}
\newtheorem{definition}{Definition}
\newtheorem{notation}{Notation}
\newtheorem{cl}{Claim}

\newtheorem{case}{Case}

\crefname{clm}{Claim}{Claims}
\crefname{cl}{Claim}{Claims}
\crefname{case}{Case}{Cases}

\title[Localization Framework]{Localization: A Framework to Generalize Extremal Graph Problems}

\author*[1]{\fnm{Rajat} \sur{Adak}\orcid{https://orcid.org/0009-0001-2723-5550}}\email{rajatadak@iisc.ac.in}

\author[1]{\fnm{L. Sunil} \sur{Chandran}}
\email{sunil@iisc.ac.in}

\affil[1]{\orgdiv{Department of Computer Science and Automation},
\orgname{Indian Institute of Science},
\orgaddress{\city{Bangalore}, \country{India}}}

\abstract{
Extremal graph theory studies the maximum or minimum number of subgraphs isomorphic to a prescribed graph under given constraints. \textit{Localization} has recently emerged as a framework that refines such problems by assigning extremal quantities locally (to vertices or edges) and then aggregating them. This perspective not only recovers classical results but also leads to sharper bounds. A classical result states that a connected planar graph with a finite girth $g$ satisfies
\[m \leq \frac{g}{g-2}(n-2)\]
Wood~\cite{wood} derived upper bounds on the number of $K_t$-cliques in graphs of bounded maximum degree, expressed in terms of both the number of vertices and the number of edges: 
\begin{align*}
    ex(n,K_t,K_{1,d+1}) \leq \frac{n}{d+1}\binom{d+1}{t} \\
mex(m,K_t,K_{1,d+1}) \leq \frac{m}{\binom{d+1}{2}}\binom{d+1}{t}
\end{align*}
More recently, Chakraborty and Chen~\cite{CHAKRABORTI2024103955} established a similar upper bound for graphs with bounded path length:  
\begin{equation*}
    mex(m,K_t,P_{r+1}) \leq \frac{m}{\binom{r}{2}}\binom{r}{t}
\end{equation*}

In this paper, we employ the localization framework to improve these bounds and provide structural characterizations of the extremal graphs attaining them.
}

\begin{document}

\maketitle

\section{Introduction}
Extremal graph theory, initiated by Turán, studies the maximum or minimum number of copies of a certain subgraph in a graph that does not contain any member of a prescribed family as a subgraph. Some of the earliest foundational results in extremal graph theory are the following classical theorems.  

\begin{theorem}\label{th:Turan}\emph{(Tur\'{a}n \cite{Turan})}
Let $G$ be a simple graph on $n$ vertices with clique number at most $r$. Then
\begin{equation*}
    |E(G)| \leq \frac{n^2(r-1)}{2r},
\end{equation*}
with equality if and only if $G$ is a regular Tur\'{a}n graph on $n$ vertices with $r$ classes.
\end{theorem}
\begin{definition}
A planar graph $G$ is called a \emph{$k$-angulation} if it admits a planar embedding in which every face (including the outer face) is a cycle of length $k$.
\end{definition}
\begin{theorem}\label{thm:girth}
Let $G$ be a simple planar connected graph on $n$ vertices with finite girth $g$, where girth is the length of the smallest cycle in the graph. Then
\[|E(G)| \leq \frac{g}{g-2}(n-2)\]
with equality if and only if $G$ is a $g$-angulation.
\end{theorem}

Let $\mathcal{F}$ be a family of graphs. A graph is considered $\mathcal{F}$-free if it contains no subgraph that is isomorphic to any graph in the family $\mathcal{F}$. Tur\'{a}n number of an $\mathcal{F}$-free graph on $n$ vertices is denoted by $ex(n,\mathcal{F})$, which counts the maximum number of edges in such a graph. The generalized Tur\'{a}n numbers $ex(n,H,\mathcal{F})$ and $mex(m,H,\mathcal{F})$ extend this framework by counting the maximum number of copies of a fixed graph $H$ that can appear in $\mathcal{F}$-free graphs on $n$ vertices or with $m$ edges, respectively. For the special case $\mathcal{F} = \{F\}$, we simplify the notation to $ex(n,H,F)$ (or $mex(m,H,F)$).
These parameters have been widely studied when $\mathcal{F}$ is a natural family of graphs, such as paths, cycles, and stars and $H$ is a clique.
\begin{notation}
Let $N(G, K_t)$ denote the number of subgraphs of $G$ isomorphic to $K_t$.
   
\end{notation}
\begin{theorem}\label{wood1}\emph{(Wood~\cite{wood})} Let $t \geq 1$ and $G$ be a graph on $n$ vertices with maximum degree, $\Delta(G) = d$. Then;
\[N(G,K_t) \leq ex(n,K_t,K_{1,d+1})\leq \frac{n}{d+1}{d+1 \choose t} = \frac{n}{t}{d\choose t-1} \]
If $t = 1$, we always get equality; otherwise, equality holds if and only if $G$ is a disjoint union of copies of $K_{d+1}$. 
\end{theorem}
\begin{corollary}\label{WoodCor2}\emph{(Wood~\cite{wood})}
    Let $t \geq 1$ and $G$ be a graph on $n$ vertices with maximum degree $\Delta(G) = d$. Then  
    \[
        N(G,K_{\geq 1}) \leq \frac{n}{d+1}(2^{d+1}-1).
    \]
\end{corollary}
\begin{theorem}\label{wood2}\emph{(Wood~\cite{wood})} Let $t \geq 2$ and $G$ be a graph with $m$ edges and maximum degree, $\Delta(G) = d$. Then;
\[N(G,K_t) \leq mex(m,K_t,K_{1,d+1})\leq \frac{m}{{d+1 \choose 2}}{d+1 \choose t} = \frac{m}{{t\choose 2}}{d-1\choose t-2}\]
If $t = 2$, we always get equality; otherwise, equality holds if and only if $G$ is a disjoint union of copies of $K_{d+1}$ with any number of isolated vertices.
\end{theorem}
\begin{corollary}\label{WoodCor1}\emph{(Wood~\cite{wood})}
    Let $t \geq 2$ and $G$ be a graph with $m$ edges and maximum degree $\Delta(G) = d$. Then  
    \[
        N(G,K_{t \geq 2}) \leq \frac{m}{{d+1 \choose 2}}(2^{d+1}-d-2).
    \]
\end{corollary}
\begin{theorem}\label{chakchen}\emph{(Chakraborti-Chen~\cite{CHAKRABORTI2024103955})} Let $t \geq 2$ and $G$ be a $P_{r+1}$-free graph on $m$ edges, then
\[N(G,K_t) \leq mex(m,K_t,P_{r+1})\leq \frac{m}{{r\choose 2}}{r\choose t}\]
If $t =2$, equality always holds; otherwise, equality holds if and only if $G$ is a disjoint union of copies of $K_r$ with any number of isolated vertices.
    
\end{theorem}

\subsection{Localization}

A classical lower bound on the independence number $\alpha(G)$ of a graph $G$ is 
\begin{equation*}
    \alpha(G) \geq \frac{n}{\Delta + 1},
\end{equation*}
where $\Delta$ denotes the maximum degree of $G$. While this inequality depends only on a global parameter, namely the maximum degree, stronger bounds can often be obtained by incorporating local structural information. In this direction, Caro~\cite{caro1979new} and Wei~\cite{wei1981lower} independently established a celebrated refinement by replacing $\Delta$ with the degree of each individual vertex. Their result asserts that
\begin{equation*}
    \alpha(G) \geq \sum_{v \in V(G)} \frac{1}{d(v) + 1},
\end{equation*}
where $d(v)$ denotes the degree of vertex $v \in V(G)$. This formulation highlights the power of localized parameters in strengthening classical extremal results, and naturally, this led to the question of whether other graph parameters, when localized, could provide comparable extensions of classical extremal theorems.

    Brada\v{c}~\cite{bradac} and Malec--Tompkins~\cite{DBLP:journals/ejc/MalecT23} recently introduced the notion of \emph{localization} as a tool to strengthen and generalize classical extremal bounds. Observe that the constraint in \Cref{th:Turan} relies on a global graph parameter—specifically, the maximum clique size. In contrast, the localization framework assigns \emph{local weights} to graph elements (such as vertices or edges), thereby transforming global restrictions into locally defined conditions.

Brada\v{c}~\cite{bradac} and Malec--Tompkins~\cite{DBLP:journals/ejc/MalecT23} defined a weight function $w : E(G) \rightarrow \mathbb{Z}$, where the value of $w(e)$ is determined by the particular extremal bound under consideration. In the case of \Cref{th:Turan}, $w(e)$ is defined as the order of the largest clique in $G$ containing the edge $e$. They proved the following result:

\begin{theorem}\label{thm:GTuran}\emph{(Brada\v{c}~\cite{bradac}, Malec--Tompkins~\cite{DBLP:journals/ejc/MalecT23})}
For a simple graph $G$ with $n$ vertices,
\begin{equation*}
    \sum_{e \in E(G)} \frac{w(e)}{w(e) - 1} \leq \frac{n^2}{2},
\end{equation*}
with equality if and only if $G$ is a regular Tur\'{a}n graph.
\end{theorem}
    
    It is easy to verify that \Cref{thm:GTuran} implies the bound of \Cref{th:Turan}, when we assume $G$ to be $K_{r+1}$-free. These localized weights allow for the improvement of the classical bound in \Cref{th:Turan} by replacing the global parameter, maximum clique size, with these edge weights.

In~\cite{zhao2025localized}, the localization framework was further applied to generalize the Erd\H{o}s--Gallai theorem for cycles~\cite{gallai1959maximal} and the corresponding path theorem was localized in~\cite{DBLP:journals/ejc/MalecT23}. Kirsch and Nir~\cite{Kirsch_2024} established a localized version of Zykov’s generalization \cite{Zykov1949-qr} of Turán’s theorem, which was subsequently further improved in~\cite{aragao2024localised}, the authors provided a localization for the Graph Maclaurin Inequalities~\cite{khadzhiivanov1977inequalities}. Moreover, the localization framework has also been extended to hypergraphs~\cite{zhao2025localizedhyper} and to spectral extremal bounds~\cite{liu2026local,kannan2025localizationspectralturantypetheorems}. This further illustrates the versatility and potential of the localization framework to extend across different domains within extremal graph theory.

While edge-based localization often produces distribution-type inequalities reminiscent of the LYM or Kraft inequalities, a more direct alternative was proposed in~\cite{adak2025vertex}, where localization is applied to vertices rather than edges.

Adak and Chandran~\cite{adak2025vertex} introduced the concept of \emph{vertex-based localization}. They provided a vertex-based localization of the Erd\H{o}s--Gallai theorems for paths and cycles. In subsequent work~\cite{adak2025turan}, the same authors established a vertex-based localization of a popular form of Turán’s theorem. For a simple graph $G$ they proved that,
\[|E(G)| \leq \frac{n}{2}\left\lfloor\sum_{v \in V(G)}\frac{c(v)-1}{c(v)}\right\rfloor\] where the weight function $c: V(G) \rightarrow \mathbb{Z}$ is defined such that $c(v)$ is the order of the largest clique containing the vertex $v$.

More recently, in \cite{adak2025vertexbasedlocalizationgeneralizedturan}, they extended this localization framework to improve the bounds of generalized Erd\H{o}s--Gallai theorems due to Luo~\cite{luo2018maximum} and in \cite{adak2025generalizedzykovstheorem} they provided a vertex-based localization of Zykov's Theorem.

In this paper, we use the localization framework to generalize \Cref{thm:girth},\ref{wood1}, \ref{wood2} and \ref{chakchen}, and characterize the extremal graphs for these generalized bounds.

\section{Main Results}
We provide a vertex-based localization of \Cref{wood1}, and generalize \Cref{thm:girth,wood2,chakchen}, by assigning suitable weights to the edges of the graph and applying the localization framework.
\begin{notation}
    Define $g:E(G)\rightarrow \mathbb{Z}$, such that $g(e)$ is the length of the smallest cycle edge $e$ is part of, and if $e$ does not participate in any cycle, set $g(e) =2$.
\end{notation}

\begin{theorem}\label{girthlocal}
    Let $G$ be a planar connected graph on $n$ vertices. Then;
    \[\sum_{e \in E(G)}\frac{g(e)-2}{g(e)} \leq n-2\]
    equality holds if and only if $G$ is either $K_2$ or is a $k$-angulation for some $k$.
\end{theorem}
\begin{theorem}\label{thm2}
    Let $t \geq 1$ and $G$ be a simple graph then, 
    \[ N(G,K_t) \leq \sum_{v \in V(G)} \frac{{d(v)+1 \choose t}}{d(v)+1} = \frac{1}{t}\sum_{v \in V(G)}{d(v) \choose t-1}\]
    Let $X = \{v \in V(G) \mid d(v) \geq t-1\}$. If $t =1$, equality always holds; otherwise, equality holds if and only if all the components of $G[X]$ are cliques.
\end{theorem}
\begin{corollary}\label{woodcor1}
    Let $t \geq 1$ and $G$ be a simple graph. Then  
    \[
        N(G,K_{t \geq 1}) \leq \sum_{v \in V(G)}\sum_{t =1}^{d(v)+1}\frac{{d(v)+1 \choose t}}{d(v) +1} 
        = \sum_{v \in V(G)}\frac{2^{d(v)+1} -1}{d(v)+1}.
    \]
\end{corollary}
The above result can be seen as a localized version of \Cref{WoodCor2}.
\begin{definition}
        A graph $G$ is called a \textit{diamond-graph} if it is isomorphic to $K_4$ minus one edge. A graph is \textit{diamond-free} if it has no diamond as an induced subgraph.
\end{definition}

\begin{notation}
    Define $w: E(G) \rightarrow \mathbb{Z}$, such that $w(e)=$ Number of common neighbors of the endpoints of $e$. In other words, $w(e)$ is the number of triangles that edge $e$ is part of in $G$.  This can be seen as the \textit{degree} of edge $e$.
\end{notation}
\begin{theorem}\label{thm1}
    Let $t\geq 2$ and $G$ be a simple graph then,
    \[N(G,K_t) \leq \sum_{e\in E(G)} \frac{{w(e) +2 \choose t}}{{w(e)+2 \choose 2}} = \frac{1}{{t\choose 2}}\sum_{e \in E(G)}{w(e)\choose t-2}\]
    Let $Y = \{e \in E(G)\mid w(e) < t-2\}$. If $t =2$, equality always holds; otherwise, equality holds if and only if $G\setminus Y$ is a diamond-free graph.
\end{theorem}
\begin{corollary}\label{woodcor2}
    Let $t \geq 2$ and $G$ be a simple graph. Then  
    \[
        N(G,K_{t \geq 2}) \leq \sum_{e\in E(G)} \sum_{t =2}^{w(e)+2}\frac{{w(e)+2\choose t}}{{w(e)+2\choose 2}} 
        = \sum_{e \in E(G)}\frac{2^{w(e)+2} - w(e) -3}{{w(e) +2\choose 2}}.
    \]
\end{corollary}
\Cref{woodcor2} can be seen as a localized version of \Cref{WoodCor1}.

\begin{notation}
    Define $p: E(G) \rightarrow \mathbb{Z}$, such that $p(e)$ is the length of the longest path in $G$ containing the edge $e$.
\end{notation}
\begin{theorem}\label{thm3}
    Let $t \geq 2$ and $G$ be a simple graph, then
    \[ N(G,K_t) \leq \sum_{e \in E(G)} \frac{{p(e) +1 \choose t}}{{p(e) +1 \choose 2}} = \frac{1}{{t\choose 2}}\sum_{e\in E(G)}{p(e)-1\choose t-2}\]
    Let $Z = \{e \in E(G)\mid p(e) < t-1\}$. If $t =2$, equality always holds; otherwise, equality holds if and only if all the components of $G \setminus Z$ are cliques.
\end{theorem}
\subsection{Recovering the classical bounds}

\noindent \textbf{From \Cref{girthlocal} to \Cref{thm:girth} (For $2$-connected graphs)}
\newline Assume that $G$ is $2$-connected and has a finite girth $g$. Thus there are no cut-edges in $G$. Therefore, $g \leq g(e)$ for all $e \in E(G)$. Note that $g >2$. Therefore, from the bound of \Cref{girthlocal} we get;
\begin{align*}
    &n-2 \geq \sum_{e\in E(G)}\frac{g(e)-2}{g(e)} \geq \sum_{e\in E(G)}\frac{g-2}{g} = \frac{g-2}{g}m\\
    &\implies  m \leq \frac{g}{g-2}(n-2)
\end{align*}
For equality, we must have equality in the bound of \Cref{girthlocal}. Thus $G$ must be a $K_2$ or a $k$-angulation. Since $G$ has finite girth $G \not\cong K_2$ and since the girth of $G$ is $g$, we get that $k = g$. Thus $G$ is $g$-angulation.

\vspace{2mm}

\noindent \textbf{From \Cref{thm2} to \Cref{wood1}}
    \newline Since $\Delta(G) = d$, therefore $d \geq d(v)$ for all $v \in V(G)$. Thus $\frac{1}{t}{d(v)\choose t-1} \leq \frac{1}{t}{d\choose t-1}$ for all $v \in V(G)$. Therefore, from the bound of \Cref{thm2} we get;
    \begin{equation*}
        N(G,K_t) \leq \frac{1}{t}\sum_{v\in V(G)}{d(v) \choose t-1} \leq \frac{1}{t}\sum_{v \in V(G)}{d \choose t-1} = \frac{n}{t}{d\choose t-1}
    \end{equation*}
For equality, both the above inequalities must be equalities. Thus $G$ must be $d$-regular graph. From \Cref{thm2}, the equality holds only if all the components of $G[X]$ are cliques. Therefore, $G$ must be a disjoint union of copies of $K_{d+1}$.

\vspace{2mm}

\noindent\textbf{From \Cref{thm1} to \Cref{wood2}}
\newline Since $\Delta(G) =d$, for all $e = \{u,v\} \in E(G)$, $w(e) \leq  \text{min}\{d(v),d(u)\}-1\leq \Delta(G) -1 = d-1$. Therefore, for all $e \in E(G)$ we have, $\frac{1}{{t\choose 2}}{w(e)\choose t-2} \leq \frac{1}{{t\choose 2}}{d-1 \choose t-2}$. Thus, from the bound of \Cref{thm1} we get;
\begin{equation*}
    N(G,K_t) \leq \frac{1}{{t\choose 2}}\sum_{e\in E(G)}{w(e)\choose t-2} \leq  \frac{m}{{t\choose 2}}{d-1 \choose t-2}
\end{equation*}
For equality, we must have $w(e) = d-1$ for all $e \in E(G)$. Since $w(e) \leq  \text{min}\{d(v),d(u)\}-1\leq d-1$, for all $x \in V(G)$ such that there exists an edge adjacent to $x$, that is, $d(x) > 0$, we get that, $d(x) = d$. Also from \cref{thm1}, $G\setminus Y$ must be a diamond-free graph. Thus $G\setminus Y$ must be a disjoint union of copies of $K_{d+1}$ and since $w(e) = d-1$ for all $e \in E(G)$ and $d(x)=d$ for all $x \in V(G\setminus Y)$. Note that $Y$ must consist of only isolated vertices. Thus $G$ is a disjoint union of copies of $K_{d+1}$ and isolated vertices.
\vspace{2mm}

 \noindent\textbf{From \Cref{thm3} to \Cref{chakchen}}
\newline Since $G$ is $P_{r+1}$-free, we get that $p(e)+1 \leq r$ for all $e \in E(G)$. Since $p(e)+1 \geq 2$ we get that, $\frac{1}{{t\choose 2}}{p(e)-1\choose t-2} \leq \frac{1}{{t\choose 2}}{r-2\choose t-2} = \frac{{r\choose t}}{{r\choose 2}}$. Therefore, from the bound of \Cref{chakchen} we get;
\begin{equation*}
    N(G,K_t) \leq \frac{1}{{t\choose 2}}\sum_{e\in E(G)}{p(e)-1\choose t-2} \leq  \sum_{e \in E(G)}\frac{{r\choose t}}{{r\choose 2}} = \frac{m}{{r\choose 2}}{r\choose t}
\end{equation*}
For equality, both the above inequalities must be equalities. From \Cref{thm3}, the equality holds only if all the components of $G\setminus Z$ are cliques. Also, $p(e)+1 = r$ for all $e \in E(G)$. Thus, the cliques are either of order $r$ or $1$.
\section{Proof of Results}
\subsection{Proof of \Cref{girthlocal}}
\begin{proof}
Since $G$ is planar and connected, from Euler's formula, we get that; 
\begin{equation}\label{euler}
    |F(G)| = 2 - |V(G)| + |E(G)|
\end{equation}
Let $F(G) = \{ f_1, f_2, \dots , f_{|F(G)|}\}$. If an edge $e \in E(G)$ is incident to a face $f_i$ for some $i \in [\,|F(G)|\,]$, then $g(e) \leq d(f_i)$, where $d(f_i)$ denotes the degree of the face $f_i$, that is the number of edges in the boundary walk of $f_i$ (thus cut-edges are counted twice). 

Observe that every edge of $G$ which is not a cut edge is incident to exactly two faces, contributing once to the degree of each of them. In contrast, a cut-edge is incident to a single face but contributes twice to its degree.

\begin{figure}[htb]
    \centering 
    \begin{tikzpicture}[
        V/.style={circle, fill=black, inner sep=1.5pt},
        E/.style={draw, thin},
        font=\small
    ]

        \node[V] (v1) at (0, 1.5) {};
        \node[V] (v2) at (0, 0) {};
        \node[V] (v3) at (1.5, 0.75) {};
        \node[V] (v4) at (3.5, 0.75) {};
        \node[V] (v5) at (5, 1.5) {};
        \node[V] (v6) at (5, 0) {};

        \path[E] (v1) edge["$e_1$" above left] (v2);
        \path[E] (v2) edge["$e_2$" below] (v3);
        \path[E] (v3) edge["$e_3$" above right] (v1);

        \draw[E] (v3) -- (v4) node[midway, above, font=\small] {$e_4$};

        \path[E] (v4) edge["$e_5$" below] (v6);
        \path[E] (v6) edge["$e_6$" below right] (v5);
        \path[E] (v5) edge["$e_7$" above left] (v4);

        \node at (0.5, 0.75) {$f_1$};
        \node at (4.5, 0.75) {$f_2$};
        \node at (2.5, 1.5) {$f_3$};

    \end{tikzpicture}
    
    \caption{$E(f_1) = \{e_1,e_2,e_3\}$, thus $d(f_1) = 3$. But $E(f_3) = \{e_1,e_2,e_4,e_5,e_6,e_7,e_4,e_3\}$, note that $e_4$ comes twice, and therefore $d(f_3) = 8$.}
    \label{fig:ex} 
\end{figure}

Let $E(f_i)$ be the multiset of edges in the boundary walk of $f_i$, where the cut-edges in the boundary of $f_i$, appear twice. Refer to \cref{fig:ex} for an example. Note that $|E(f_i)| = d(f_i)$. Thus, we get;
\begin{equation}\label{maingirth}
    2 \sum_{e \in E(G)}\frac{1}{g(e)} \geq \sum_{i =1}^{|F(G)|}\sum_{e \in E(f_i)}\frac{1}{d(f_i)} = \sum_{i=1}^{|F(G)|}d(f_i)\frac{1}{d(f_i)} = |F(G)|
\end{equation}
 Thus from \cref{euler,maingirth} we get that;
 \begin{align*}
     &2\sum_{e\in E(G)}\frac{1}{g(e)} \geq 2 - |V(G)| + |E(G)|\\
     \implies &n -2 \geq \sum_{e \in E(G)}\left(1-\frac{2}{g(e)}\right) = \sum_{e \in E(G)}\left(\frac{g(e)-2}{g(e)}\right)
 \end{align*}

For the extremal analysis, we first consider the case when $G$ is either $K_{2}$ or a $k$-angulation. 
If $G \cong K_{2}$, then $E(G) = \{e\}$ and $g(e)=2$. In this case, both sides of the bound evaluate to $0$. 
If $G$ is a $k$-angulation, it follows that $g(e) = k$ for all $e \in E(G)$ and $d(f_i) = k$ for all $f_i \in F(G)$. 
Hence, equality holds in \cref{maingirth}, and consequently, in the bound of \Cref{girthlocal}.  

For the converse, assume equality holds in the bound of \Cref{girthlocal} with $G \not\cong K_{2}$. 
This implies equality in \cref{maingirth}, and therefore $d(f_i) = g(e)$ for every $f_i \in F(G)$ and every $e \in E(f_i)$. 
Since $g(e) = 2$ if and only if $e$ is a cut-edge, while $d(f_i) > 2$ for all $f_i \in F(G)$ (as $G \not\cong K_{2}$), it follows that $G$ contains no cut-edges and therefore, each face in $G$ is bounded by a cycle.

Two faces are \emph{adjacent} if they share at least one edge. In the dual graph, this corresponds to the adjacency of the associated vertices. Similarly, two faces are \emph{connected} if the corresponding vertices in the dual graph are connected. 

If $f_i$ and $f_j$ are adjacent and share an edge $e$, then we have $d(f_i) = d(f_j) = g(e)$. Since the dual graph of any planar graph is connected, every face is connected to every other face. By propagating this equality of $d(\cdot)$ along paths in the dual, it follows that all faces share the same degree. Hence $d(f)$ is constant for all $f \in F(G)$.

Since every face is bounded by a cycle, the existence of a face $f$ whose boundary cycle has length strictly less than $d(f)$ would imply the presence of an edge $e$ on the boundary of $f$ with $g(e) < d(f)$. This is a contradiction. Hence, each face must be bounded by a simple cycle of length exactly $d(f)$. Consequently, $G$ is a $k$-angulation with $k = d(f)$ for any $f \in F(G)$.
\end{proof}

\subsection{Proof of \Cref{thm2}}\label{vertexdegreeproof}
\begin{notation}
    If $v \in V(G)$, then $N(G,K_t,v)$ denotes the number of copies of $K_t$ in $G$ containing the vertex $v$.
\end{notation}
Recall that the set $X$ consists of all vertices in $G$ whose degrees are at least $t-1$.
\begin{lemma}\label{lem1}
\Cref{thm2} holds for $G$ is and only if it holds for $G\setminus X$.
\end{lemma}
\begin{proof}This follows from the fact that vertices in the complement set $X^c=V(G) \setminus X$ do not contribute to either side of the inequality in \Cref{thm2}. Note that $N(G,K_t) = N(G[X],K_t)$, since any clique of order at least $t$ in $G$, that is not entirely in $G[X]$ has to contain a vertex $v \in X^c$; but then $d(v) \geq t-1$, contradiction. On the other hand, since for $v \in X^c$, we have $d(v) < t-1$, we get ${d(v) \choose t-1} = 0$. 
 \end{proof}
Therefore, without loss of generality, we may assume that $G$ contains only vertices of degree at least $t-1$. Now we are ready to prove \Cref{thm2}.
    \begin{proof}
        For $t=1$, the bound and the equality case can both be verified trivially; thus, we assume $t \geq 2$. Since we assumed $d(v) \geq t-1$ for all $v \in V(G)$, we get that $G$ does not contain any isolated vertex; therefore, $\delta(G) > 0$.
        \newline Proof by induction on the number of vertices.
        \vspace{2mm}

        \noindent \textit{Base Case:} For $n=2$, let $V(G)=\{u,v\}$ and $E(G) = \{\{u,v\}\}$ (since $\delta(G) >0)$. For $t>2$, the bound is trivially true since both sides are zero. For $t =2$, $N(G,K_t)=1$ and the right-hand side also sums up to $1$.
        \vspace{1mm}

        \noindent\textit{Induction Hypothesis:} Suppose the claim in \Cref{thm2} is true for all graphs with less than $|V(G)|$ vertices.
        \vspace{1mm}

        \noindent\textit{Induction Step:} Let $x \in V(G)$ such that $0<d(x) = \delta(G)$. Now define $G' = G\setminus \{x\}$. Recall that $N(G,K_t,x)$ denotes the number of copies of $K_t$ in $G$ containing the vertex $x \in V(G)$. Clearly;
    \begin{equation}\label{eqn4}
         N(G,K_t) = N(G',K_t) + N(G,K_t,x)
    \end{equation}
    Note that for every $K_t$ containing $x$, there exists a unique copy of $K_{t-1}$ in $G[N(x)]$. Therefore;
\begin{equation}\label{eqn5}
  N(G,K_t,x) = N(G[N(x)],K_{t-1}) \leq {d(x) \choose t-1}
\end{equation}
Let $d_{G'}(v)$ be the degree of vertex $v$ in $G'$. Now note that $d_{G'}(v) = d(v)$ for all $v\in V(G') \setminus N(x)$, where $N(x)$ denotes the set of vertices adjacent to $x$ in $G$. Clearly, $d_{G'}(v) = d(v) -1$ for all $v \in N(x)$. From the induction hypothesis, we know that;
\begin{equation}\label{eqn6}
    N(G',K_t) \leq \frac{1}{t}\sum_{v \in V(G')}{d_{G'}(v) \choose t-1} = \frac{1}{t}\sum_{v \in V(G')\setminus N(x)}{d(v) \choose t-1} + \frac{1}{t}\sum_{v \in N(x)} {d(v)-1 \choose t-1}
\end{equation}
Thus from \cref{eqn4,eqn5,eqn6} we get that;
\begin{equation}\label{eqn7}
    N(G,K_t) \leq \frac{1}{t}\sum_{v \in V(G')\setminus N(x)}{d(v) \choose t-1} + \frac{1}{t}\sum_{v \in N(x)} {d(v)-1 \choose t-1} +{d(x) \choose t-1}
\end{equation}
Note that;
\begin{equation}\label{eqn8}
    \frac{1}{t}{d(v)-1 \choose t-1} = \frac{(d(v)-t+1)}{td(v)}{d(v) \choose t-1} = \frac{1}{t}{d(v) \choose t-1} - \frac{1}{d(v)}{d(v) \choose t-1}+\frac{1}{td(v)}{d(v) \choose t-1}
\end{equation}
Thus from \cref{eqn7,eqn8} we get;
\begin{align}\label{eqn9}
    N(G,K_t) &\leq {d(x) \choose t-1}+ \frac{1}{t}\sum_{v \in V(G')\setminus N(x)}{d(v) \choose t-1} + \sum_{v \in N(x)} \left(\frac{1}{t}{d(v) \choose t-1} - \frac{1}{d(v)}{d(v) \choose t-1}+\frac{1}{td(v)}{d(v) \choose t-1}\right)\nonumber\\
    &= {d(x) \choose t-1} + \frac{1}{t}\sum_{v\in V(G')}{d(v) \choose t-1}  +  \sum_{v \in N(x)}\left(\frac{1}{td(v)}{d(v) \choose t-1} - \frac{1}{d(v)}{d(v) \choose t-1}\right) 
\end{align}
\begin{claim}\label{claim1} For $d(v) \geq d(x) >0$ and $t \geq 2$ we have;
    \[\sum_{v \in N(x)}\left(\frac{1}{td(v)}{d(v) \choose t-1} - \frac{1}{d(v)}{d(v) \choose t-1}\right) + {d(x) \choose t-1} \leq \frac{1}{t}{d(x) \choose t-1}\]
    \begin{proof}
        Note that $\delta(G) =d(x) \leq d(v)$ for all $v \in V(G)$. Therefore, we get $\frac{1}{t-1}{d(v)-1\choose t-2} \geq \frac{1}{t-1}{d(x)-1\choose t-2}\implies \frac{1}{d(v)}{d(v) \choose t-1} \geq \frac{1}{d(x)}{d(x) \choose t-1}$. Thus we get;
        \begin{align}\label{eqn10}
            & {d(x) \choose t-1} = \sum_{v \in N(x)}\frac{1}{d(x)}{d(x) \choose t-1} \leq \sum_{v \in N(x)}\frac{1}{d(v)}{d(v)\choose t-1} \\
            \implies & \frac{t-1}{t}{d(x) \choose t-1} \leq \sum_{v \in N(x)} \frac{t-1}{td(v)}{d(v)\choose t-1}\nonumber \\
            \implies & {d(x) \choose t-1} - \frac{1}{t}{d(x) \choose t-1} \leq \sum_{v \in N(x)}\left(-\frac{1}{td(v)}{d(v) \choose t-1} + \frac{1}{d(v)}{d(v) \choose t-1}\right)\nonumber\\
            \implies & \sum_{v \in N(x)}\left(\frac{1}{td(v)}{d(v) \choose t-1} - \frac{1}{d(v)}{d(v) \choose t-1}\right) + {d(x) \choose t-1} \leq \frac{1}{t}{d(x) \choose t-1}\nonumber
        \end{align}
  \end{proof}
\end{claim}
From \cref{eqn9,claim1} we get that;
\[N(G,K_t) \leq \frac{1}{t}{d(x) \choose t-1} +\frac{1}{t}\sum_{v \in V(G')}{d(v) \choose t-1} = \frac{1}{t}\sum_{v \in V(G)}{d(v) \choose t-1}\]
Thus, we get the required bound. Now we will characterize the extremal graphs for \Cref{thm2}.

First, suppose that all the components of $G$ are cliques. Let $C$ be one such component with order $k$. Clearly $d(v) = k-1$ for all $v \in V(C)$. And $C$ contains ${k \choose t}$ copies of $K_t$. The contribution by the vertices of $C$ to the left side of the bound in \Cref{thm2} is, 
\[\frac{1}{t}\sum_{v \in V(C)}{k-1\choose t-1} = \frac{k}{t}{k-1\choose t-1} = {k\choose t}\]
Thus, we get equality in the bound.

Now, assume that we have equality in the bounds of \Cref{thm2}. Thus, we must have equality in \cref{eqn5,eqn6,eqn10}. From equality in \cref{eqn5}, we get that $G[N[x]]$ is a clique, where $N[x]$ denotes the closed neighborhood of $x$. Using the induction hypothesis and the equality in \cref{eqn6}, we get that all the components of $G'$ are cliques. Assume that $C$ is a connected component of $G$, then it is enough to show that $C$ is a clique. 

If $C$ does not contain $x$, clearly $C$ does not contain any vertex from $N(x)$. Thus, $C$ is a component in $G'$ as well, and therefore, by induction hypothesis, $C$ is a clique. Now suppose $C$ contains the vertex $x$, therefore $C$ contains all the vertices of $N(x)$. Since $N[x]$ induces a clique in $G$, $C\setminus \{x\}$ must be connected and therefore, $C\setminus \{x\}$ must be a component of $G'$, and thus is a clique. Suppose there exists $y \in V(C) \setminus N[x]$. Clearly $y$ is adjacent to $v \in N(x)$ but not adjacent to $x$. Therefore, $d(v) > d(x)$. But from equality in the bound of \cref{eqn10} we get that, $d(x) = d(v)$ for all $v \in N(x)$. We get a contradiction. Therefore, $V(C) = N[x]$, which induces a clique.    \end{proof}
\subsection{Proof of \Cref{thm1}}
 Recall that $Y = \{e\in E(G) \mid w(e) +2 < t\}$
\begin{lemma}\label{lem2}
\Cref{thm1} holds for $G$ if and only if it holds for $G\setminus Y$
\end{lemma} 
   \begin{proof}
This follows from the fact that edges in the set $Y$ do not contribute to either side of the inequality in \Cref{thm1}. Note that for any edge $e \in Y$, the endpoints of $e$ have fewer than $t-2$ common neighbors in $G$. Hence, no copy of $K_t$ in $G$ can contain $e$. On the other hand, since $w(e) < t-2$, we have ${w(e) \choose t-2} = 0$. 
\end{proof}
Therefore, without loss of generality, we may assume that $G$ contains only edges with weights at least $t-2$. Now coming to the proof of \Cref{thm1}.
    \begin{proof}
    If $t =2$, clearly the bound holds with equality; thus, assume $t \geq 3$. Let $e \in E(G)$, if $e = \{u,v\}$ where $u,v \in V(G)$. If $A$ is a $K_t$ containing the edge $e$, then $u,v \in V(A)$ and $V(A)\setminus\{u,v\} \subseteq N(u) \cap N(v)$. Note that $|N(u)\cap N(v)| = w(e)$. Thus, the number of copies of $K_t$ containing the edge $e$ is at most ${w(e) \choose t-2}$. If we add up the contribution of each edge, that is the number of copies of $K_t$, containing an edge, then each copy is counted ${t\choose 2}$ times, since it is counted once for each edge participating in that $K_t$ copy. Thus, we get;
    \begin{equation}\label{e1}
        N(G,K_t) \leq \frac{1}{{t\choose 2}}\sum_{e\in E(G)}{w(e) \choose t-2}
    \end{equation}
Thus, we get the required bound. Now we will characterize the extremal graphs for \Cref{thm1}.

    First, assume we have equality in \Cref{thm1}. Therefore, we must have the number of copies of $K_t$ containing the edge $e$ to be exactly ${w(e) \choose t-2}$ for all $e \in E(G)$. Thus, the vertices in $N(u)\cap N(v)$ must induce a clique where $u,v \in V(G)$ are the endpoints of $e$. Now suppose $G$ contains an induced diamond as shown in \cref{diamond}. Note that $\{a,b\} \subseteq N(c)\cap N(d)$, but $a$ and $b$ are not adjacent. Thus, we get a contradiction; therefore, $G$ is diamond-free.

    \begin{figure}[H]
        \centering
 \begin{tikzpicture}[scale=2, every node/.style={circle, draw, minimum size=8pt, inner sep=1pt}]
    \node (a) at (0,0) {a};
    \node (b) at (1,0) {b};
    \node (c) at (0.5,0.5) {c}; 
    \node (d) at (0.5,-0.5) {d}; 

    \draw (a)--(c);
    \draw (b)--(c);
    \draw (a)--(d);
    \draw (b)--(d);
    \draw (c)--(d);
\end{tikzpicture}
        \caption{Diamond in the extremal graph $G$}
        \label{diamond}
    \end{figure}

    Now suppose $G$ is diamond-free. Suppose there exists an edge $e =\{x,y\} \in E(G)$ such that $N(x)\cap N(y)$ does not induce a clique in $G$. Since $N(x)\cap N(y)$ is non-empty, there exists $u,v \in N(x) \cap N(y)$ such that $u$ and $v$ are not adjacent. Therefore, we get a diamond induced by $u,v,x$, and $y$. Thus, the common neighbors of the endpoints of any edge in $G$ must induce a clique. Therefore, the number of copies of $K_t$ containing some $e \in E(G) = {w(e) \choose t-2}$. Therefore, we get $G$ must be extremal for \Cref{thm1}.
     \end{proof}

\subsection{Proof of \Cref{thm2} using \Cref{thm1}} While an independent proof of \Cref{thm2} was presented in \cref{vertexdegreeproof}, we now offer an alternative derivation using \Cref{thm1}. This does not diminish the value of the original argument, as the degree function $d(v)$ for $v \in V(G)$ is a more natural and widely adopted parameter compared to the \textit{edge-degree} parameter, $w(e)$ for $e \in E(G)$.
 \begin{proof}
        For $ t=1$, it is easy to verify that the bound holds with equality; thus, we will assume $ t\geq 2$. Let $e \in E(G)$ such that $\{u,v\} = e$ where $u,v \in V(G)$. Clearly $w(e) \leq d(u) -1$ and $w(e) \leq d(v) -1$. Note that the number of edges in $E(G)$ which have a vertex $x \in V(G)$ as one endpoint is $d(x)$. Thus from \Cref{thm1} we get;
        \begin{align}
            N(G,K_t) \leq & \frac{1}{{t \choose 2}}\sum_{e \in E(G)}{w(e) \choose t-2} = \frac{2}{t(t-1)}\sum_{e \in E(G)}{w(e) \choose t-2}\label{eqn12}\\
            \leq & \frac{1}{t(t-1)}\sum_{\{u,v\} \in E(G)}\left({d(u) -1 \choose t-2}  + {d(v)-1 \choose t-2}\right)\label{eqn130}\\
            = & \frac{1}{t(t-1)}\sum_{v \in V(G)}d(v){d(v)-1\choose t-2}\nonumber\\
            = & \frac{1}{t}\sum_{v \in V(G)}{d(v) \choose t-1}\nonumber
        \end{align}
        Thus, we obtain the desired inequality. Observe that equality holds if and only if both inequalities in \cref{eqn12,eqn130} are tight. Specifically, equality in \cref{eqn12} holds if and only if $G$ is diamond-free, while equality in \cref{eqn130} holds if and only if for every edge $e = \{u, v\} \in E(G)$, we have $w(e) = d(u) - 1 = d(v) - 1$. 

        Note that if $G$ is a disjoint union of cliques, it is straightforward to verify that the above inequality holds with equality.  

        For the converse, assume equality holds in \cref{eqn12,eqn130}. Let $C$ be a connected component of $G$; it suffices to show that $C$ is a clique. Let $K$ be a maximum clique in $C$. Clearly, $|V(K)| > 1$. Suppose $V(C) \setminus V(K) \neq \emptyset$. Then there exists some $v \in V(K)$ adjacent to a vertex $x \in V(C) \setminus V(K)$. Choose $u \in V(K)$ with $u \neq v$. Since $w(\{u,v\}) = d(v)-1 = d(u)-1$, we have $N[u] = N[v]$. In particular, $x$ is adjacent to $u$.  

        Thus, $x,u,v$ form a triangle, and hence $|V(K)| \geq 3$. Consequently, there exists some $z \in V(K)$ not adjacent to $x$; otherwise, $K \cup \{x\}$ would induce a larger clique in $C$, contradicting the maximality of $K$. But then $x,u,v,z$ induce a diamond in $C$, and hence in $G$, contradicting the equality condition in \cref{eqn12}. Therefore, $V(C) \setminus V(K) = \emptyset$, and $C$ must be a clique.
     \end{proof}
\subsection{Proof of \Cref{thm3}}
 Recall that $Z= \{e \in E(G) \mid p(e) +1 < t\}$.
\begin{lemma}\label{lem3}
\Cref{thm3} holds for $G$ if and only if it holds for $G\setminus Z$.
\end{lemma}
\begin{proof}
This follows from the fact that edges in the set $Z$ do not contribute to either side of the inequality in \Cref{thm3}. Note that for any edge $e \in Z$, the longest path containing $e$ has length less than $t-1$. Hence, no copy of $K_t$ in $G$ can contain $e$. On the other hand, since $p(e) < t-1$, we have $p(e)-1 < t-2$, and therefore ${p(e)-1 \choose t-2} = 0$. 
\end{proof}
Therefore, without loss of generality, we may assume $p(e)\geq t-1$ for all $e \in E(G)$.
\begin{proof}
    If $t = 2$, clearly the bound holds with equality; thus, assume $t \geq 3$.
    \newline Proof by induction on the number of vertices of the graph.
    \vspace{2mm}

    \noindent\textit{Base Case:} For $n =1$, then there are no edges and the claim is vacuously true.
    \newline For $n=2$ and $|E(G)| =1$ (For $|E(G)| = 0$, the claim is again vacuously true), say $E(G) =\{e\}$. Note that $p(e) =1$. Since $t>2$, clearly both sides of the bound are zero. For $n = 3$, let $V(G) = \{u,v,w\}$ and $ \{\{u,v\},\{v,w\}\} \subseteq E(G)$ (since $p(e) \geq t-1 \geq 2)$. If $t >3$, the bound is trivially true as both sides are zero. If $t = 3$, $N(G,K_t) \leq 1$, with equality if and only if $G \cong K_3$. Clearly, the right-hand side sums up to $1$ when $t =3$.
    \vspace{2mm}

    \noindent\textit{Induction Hypothesis:} Suppose the claim in \Cref{thm3} is true for all graphs with less that $|V(G)|$ vertices.
     \vspace{2mm}

    \noindent\textit{Induction Step:} If $G$ is disconnected, by the induction hypothesis, the statement of the theorem is true for each connected component and thus for $G$. Therefore, assume $G$ to be connected. Let $k =max\{ p(e) \mid e \in E(G)\}$ and $P$ be a $k$-length path in $G$. Let $P = v_0,v_1,\dots,v_k$. Now consider the following two cases;
     \begin{case}\label{case1}
        There exists a $(k+1)$-length cycle in $G$.
        \newline Without loss of generality assume the endpoints of $P$, $v_0$ and $v_k$, are adjacent thus forming a $(k+1)$-length cycle.
        \begin{cl}\label{cl3.4}
     $V(P) = V(G)$.
 \end{cl}
\begin{proof}
          We know that $v_0$ is adjacent to $v_k$. Suppose $V(G) \setminus V(P) \neq \emptyset$. Since $G$ is connected, there exists a vertex  $v \in V(G)\setminus V(P)$, adjacent to some vertex in $P$. Since the vertices of $P$ form a cycle, without loss of generality, we can assume $v$ is adjacent to $v_0$. Now consider the path;
         \begin{equation*}
             Q = v,v_0,v_1,\dots,v_{k}
         \end{equation*}
Thus $Q$ is a path in $G$ with length $k +1$. Thus, we arrive at a contradiction to the assumption that the length of the longest path in $G$ is $k$.  \end{proof}
 
 From \cref{cl3.4} we get that $p(e) = |V(G)|-1=n-1$, for all $e \in E(G)$. Therefore;
 \begin{equation}\label{eqn13}
     \frac{1}{{t\choose 2}}\sum_{e \in E(G)}{{p(e)-1\choose t-2}} = \frac{1}{{t\choose 2}}\sum_{e \in E(G)}{n-2\choose t-2} 
 \end{equation}
Recall from the bound of \Cref{thm1};
\begin{equation*}
    N(G,K_t) \leq \frac{1}{{t\choose 2}}\sum_{e \in E(G)}{w(e) \choose t-2}
\end{equation*}
Where $w(e)$ is the number of common neighbors of the endpoints of $e \in E(G)$. Clearly $w(e) \leq n-2$ for all $e \in E(G)$. Thus we get, $\frac{1}{{t\choose 2}}{w(e) \choose t-2} \leq  \frac{1}{{t\choose 2}}{n-2\choose t-2}$. Therefore from \cref{eqn13} we get that;
\begin{equation}\label{eqn14}
    N(G,K_t) \leq \frac{1}{{t\choose 2}}\sum_{e \in E(G)}{w(e) \choose t-2} \leq \sum_{e \in E(G)} \frac{1}{{t\choose 2}}{n-2\choose t-2} =  \frac{1}{{t\choose 2}}\sum_{e \in E(G)}{{p(e)-1\choose t-2}}
\end{equation}
Clearly, for equality $w(e) = n-2$ for all $e \in E(G)$. Thus $G$ must be a clique.
\end{case}
\begin{case}\label{case2}
        There does not exist any $(k+1)$-length cycle in $G$.
        \newline From the assumption of the case, $v_0$ and $v_k$ are non-adjacent. Without loss of generality assume $d(v_0) \geq d(v_k)$.
\begin{cl}\label{light}
    $d(v_k) \leq \frac{k}{2}$
    \begin{proof}
        Suppose not, that is $d(v_k) \geq \frac{k+1}{2}$. Thus we get $d(v_0) \geq \frac{k+1}{2}$. Since $P$ is a longest path in $G$ and $v_0$ is not adjacent to $v_{k}$, we get, $N(v_0), N(v_k) \subseteq \{v_1,v_2,\dots, v_{k-1}\}$. Let $N = N(v_k) \setminus \{v_{k-1}\}$ and $N^+ = \{v_{j+1} \mid v_j \in N\}$. Clearly, $|N| = |N^+| \geq \frac{k+1}{2}-1 = \frac{k-1}{2}$. Note that $v_0, v_k \notin N^+$. Therefore, $|\{v_1,v_2,\dots,v_{k-1}\} \setminus N^+| \leq \frac{k-1}{2}$. Since $|N(v_0)| \geq \frac{k+1}{2}$, by pigeonhole principle, we get that, there exists $v_i \in N^+$ such that $v_i \in N(v_0)$ and $v_{i-1} \in N(v_k)$. Consider the cycle; $$K = v_0, v_1, \dots, v_{i-1},v_k,v_{k-1},\dots, v_i,v_0$$ Clearly $K$ is a $(k+1)$-length cycle in $G$. Thus, we get a contradiction.
     \end{proof}
    \end{cl}
   \noindent Define $G' = G\setminus \{v_k\}$. Now, using the induction hypothesis, we get that;
    \begin{equation*}
        N(G',K_t) \leq   \frac{1}{{t\choose 2}}\sum_{e \in E(G')}{p_{G'}(e)-1\choose t-2}
    \end{equation*}
    where for $e \in E(G')$, $p_{G'}(e)$ is the weight of the edge $e$ restricted to $G'$. Clearly $p_{G'}(e) \leq p(e)$ for all $e \in E(G')$. Thus we get that;
    \begin{equation*}
        N(G',K_t) \leq   \frac{1}{{t\choose 2}}\sum_{e \in E(G')}{p_{G'}(e)-1\choose t-2} \leq \frac{1}{{t\choose 2}}\sum_{e \in E(G')}{p(e)-1\choose t-2}
    \end{equation*}
    Note that the number of copies of $K_t$ containing the vertex $v_k$ is at most ${d(v_k)\choose t-1}$. Therefore,
    \begin{equation}\label{eqn15}
        N(G,K_t) \leq {d(v_k)\choose t-1} +\frac{1}{{t\choose 2}}\sum_{e \in E(G')}{p(e)-1\choose t-2}
    \end{equation}

    \begin{cl}\label{claim3.6}
        $p(e) = k$, for all $e \in E(G)\setminus E(G')$.
    \end{cl}
    \begin{proof}
            Recall that, since $P$ is the longest path in $G$, $N(v_k) \subseteq \{v_1,v_2,\dots,v_{k-1}\}$. Define $e_i = \{v_i,v_k\} \in E(G)\setminus E(G')$ where $v_i \in N(v_k)$. Thus $E(G)\setminus E(G') = \{e_i \mid v_i \in N(v_k)\}$. Define $P_i$ to be the path formed by removing the edge $\{v_i,v_{i+1}\}$ from $P$ and adding edge $e_i$. Clearly, $P_i$ is also a $k$-length path. Therefore, $p(e) = k$ for all $e \in E(G)\setminus E(G')$.
         \end{proof}
    \begin{cl}\label{claim3.7} For $t \geq 3$ we have;
        \begin{equation*}
            {d(v_k) \choose t-1}< \frac{1}{{t\choose 2}}\sum_{e \in E(G) \setminus E(G')}{{p(e)-1\choose t-2}} 
        \end{equation*}
    \end{cl}
\begin{proof}
        If ${d(v_k)\choose t-1} = 0$, the bound is trivially true, since $p(e) \geq t-1$ for all $ e\in E(G)$. Thus we assume ${d(v_k) \choose t-1} > 0$. From \cref{claim3.6} we have that;
        \begin{equation*}
            \frac{1}{{t\choose 2}}\sum_{e \in E(G) \setminus E(G')}{{p(e)-1\choose t-2}} = \frac{d(v_k)}{{t\choose 2}}{k-1 \choose t-2}
        \end{equation*}
            Recall that from \cref{light}, $d(v_k) \leq \frac{k}{2}$. Therefore, we get;
            \begin{equation*}
               \frac{d(v_k)}{{t\choose 2}}{k-1\choose t-2} \geq \frac{d(v_k)}{{t\choose 2}}{2d(v_k)-1\choose t-2} = \frac{2d(v_k)}{t(t-1)}{2d(v_k)-1\choose t-2} = \frac{1}{t}{2d(v_k)\choose t-1}
            \end{equation*}
            Since we assumed $t \geq 3$, we get that;
            \begin{equation*}
                \frac{\frac{d(v_k)}{{t\choose 2}}{k-1\choose t-2}}{{d(v_k)\choose t-1}} \geq \frac{\frac{1}{t}{2d(v_k)\choose t-1}}{{d(v_k)\choose t-1}} = \frac{2d(v_k)(2d(v_k)-1)\dots(2d(v_k)-t+2)}{td(v_k)(d(v_k)-1)\dots(d(v_k)-t+2)}\geq \frac{2^{t-1}}{t}> 1
            \end{equation*}
         \end{proof}
    
    From \cref{eqn15,claim3.7} we get that;
    \begin{equation*}
        N(G,K_t) < \frac{1}{{t\choose 2}}\sum_{e \in E(G')}{p(e)-1 \choose t-2} + \frac{1}{{t\choose 2}}\sum_{e \in E(G) \setminus E(G')}{{p(e)-1\choose t-2}} = \frac{1}{{t\choose 2}}\sum_{e \in E(G)}{{p(e)-1\choose t-2}} 
    \end{equation*}
    Clearly, we can not get equality in the bounds of \Cref{thm3} for \cref{case2}
    \end{case}
   
Equality in \Cref{thm3} is possible only for \cref{case1}. Thus, we must have equality in \cref{eqn14} and therefore $w(e) =n-2$ for all $e \in E(G)$; thus, $G$ must be a clique. 
  \end{proof}

 \section{Concluding Remarks}\label{end}
In this paper, we developed localized versions of classical extremal bounds, focusing on the number of edges in planar graphs and clique bounds for graphs with bounded maximum degree and bounded path length. By introducing specifically designed local parameters, our results refine and extend the classical bound for planar graphs and improve the bounds of Wood~\cite{wood} and Chakraborty--Chen~\cite{CHAKRABORTI2024103955}. We also characterized the extremal graphs that attain these bounds. These results demonstrate the strength of the localization framework: by shifting the analysis from global to local parameters, it not only recovers known extremal results but also provides sharper bounds and new structural insights, underscoring its potential as a powerful tool for generalizing a wide range of extremal problems.

\section*{Statements and Declarations}

\noindent\textbf{Funding} No funds, grants, or other support were received.

\noindent\textbf{Data Availability} No data was used for the research described in this article.

\noindent\textbf{Competing Interests} The authors declare that they have no known competing financial interests
or personal relationships that could have appeared to influence the work
reported in this paper.

\bibliography{references}

@article{wood,
  title={On the maximum number of cliques in a graph},
  author={Wood, David R},
  journal={Graphs and Combinatorics},
  volume={23},
  number={3},
  pages={337--352},
  year={2007},
  publisher={Springer},
  doi={10.1007/s00373-007-0738-8}
}

@article{CHAKRABORTI2024103955,
  title = {{Exact results on generalized Erdős-Gallai problems}},
  journal = {European Journal of Combinatorics},
  volume = {120},
  pages = {103955},
  year = {2024},
  issn = {0195-6698},
  doi = {https://doi.org/10.1016/j.ejc.2024.103955},
  author = {Debsoumya Chakraborti and Da Qi Chen}
}

@article{luo2018maximum,
  title={The maximum number of cliques in graphs without long cycles},
  author={Luo, Ruth},
  journal={Journal of Combinatorial Theory, Series B},
  volume={128},
  pages={219--226},
  year={2018},
  publisher={Elsevier},
  doi={10.1016/j.jctb.2017.08.005}
}

@article{DBLP:journals/ejc/MalecT23,
  author       = {David Malec and Casey Tompkins},
  title        = {Localized versions of extremal problems},
  journal      = {Eur. J. Comb.},
  volume       = {112},
  pages        = {103715},
  year         = {2023},
  doi          = {10.1016/J.EJC.2023.103715},
  bibsource    = {dblp computer science bibliography, https://dblp.org}
}

@article{bradac,
  title={A generalization of {Turán's} theorem},
  author={Bradač, Domagoj},
  journal={arXiv preprint arXiv:2205.08923},
  year={2022},
  doi={10.48550/arXiv.2205.08923}
}

@article{gallai1959maximal,
  title={On maximal paths and circuits of graphs},
  author={Erd{\H{o}}s,  P{\'a}l and Gallai, Tibor},
  journal={Acta Math. Acad. Sci. Hungar. v10},
  pages={337--356},
  year={1959},
  doi={10.1007/bf02024498}
}

@article{zhao2025localized,
  title={A localized approach for Tur{\'a}n number of long cycles},
  author={Zhao, Kai and Zhang, Xiao-Dong},
  journal={Journal of Graph Theory},
  volume={108},
  number={3},
  pages={582--607},
  year={2025},
  publisher={Wiley Online Library},
  doi={10.1002/jgt.23191}
}

@article{Turan,
  author    = {Tur\'{a}n, Paul},
  title     = {On an extremal problem in graph theory},
  journal   = {Matematikai és Fizikai Lapok},
  volume    = {48},
  pages     = {436--452},
  year      = {1941},
  language  = {Hungarian}
}

@techreport{caro1979new,
  title={New results on the independence number},
  author={Caro, Yair},
  year={1979},
  institution={Technical Report, Tel-Aviv University}
}

@misc{wei1981lower,
  title={A lower bound on the stability number of a simple graph},
  author={Wei, Victor K},
  year={1981},
  publisher={Bell Laboratories Technical Memorandum Murray Hill, NJ, USA}
}

@article{Kirsch_2024,
   title={{A Localized Approach to Generalized Turán Problems}},
   volume={31},
   ISSN={1077-8926},
   DOI={10.37236/12132},
   number={3},
   journal={The Electronic Journal of Combinatorics},
   publisher={The Electronic Journal of Combinatorics},
   author={Kirsch, Rachel and Nir, JD},
   year={2024},
   month=sep 
}

@article{aragao2024localised,
  title={{Localised graph Maclaurin inequalities}},
  author={Arag{\~a}o, Lucas and Souza, Victor},
  journal={Annals of Combinatorics},
  volume={28},
  number={3},
  pages={1021--1033},
  year={2024},
  publisher={Springer},
  doi={10.1007/s00026-023-00672-0}
}

@article{khadzhiivanov1977inequalities,
  title={Inequalities for graphs},
  author={Khadzhiivanov, N},
  journal={CR Acad. Sci. Bul},
  volume={30},
  number={6},
  pages={793},
  year={1977}
}

@article{Zykov1949-qr,
      author={Zykov, A A},
       title={{On some properties of linear complexes}},
        year={1949},
     journal={Matematicheskii Sbornik, Novaya Seriy}
}

@article{adak2025vertex,
  title={{Vertex-Based Localization of Erd\H{o}s-Gallai Theorems for Paths and Cycles}},
  author={Adak, Rajat and Chandran, L Sunil},
  journal={arXiv preprint arXiv:2504.01501},
  year={2025},
  doi={10.48550/arXiv.2504.01501}
}

@article{zhao2025localizedhyper,
  title={{Localized version of hypergraph Erd{\H{o}}s-Gallai Theorem}},
  author={Zhao, Kai and Zhang, Xiao-Dong},
  journal={Discrete Mathematics},
  volume={348},
  number={1},
  pages={114293},
  year={2025},
  publisher={Elsevier},
  doi={10.1016/j.disc.2024.114293}
}

@article{liu2026local,
  title={Local properties of the spectral radius and Perron vector in graphs},
  author={Liu, Lele and Ning, Bo},
  journal={Journal of Combinatorial Theory, Series B},
  volume={176},
  pages={241--253},
  year={2026},
  publisher={Elsevier},
  doi={10.1016/j.jctb.2025.09.001}
}

@article{adak2025turan,
  title={{Vertex-Based Localization of Tur\'{a}n's Theorem}},
  author={Adak, Rajat and Chandran, L Sunil},
  journal={arXiv preprint arXiv:2504.02806},
  year={2025},
  doi={10.48550/arXiv.2504.02806}
}

@article{kannan2025localizationspectralturantypetheorems,
      title={{Localization of spectral Tur\'{a}n-type theorems}}, 
      author={M. Rajesh Kannan and Hitesh Kumar and Shivaramakrishna Pragada},
      year={2025},
      journal={arXiv preprint arXiv:2512.01409},
      doi={10.48550/arXiv.2512.01409}
}

@article{adak2025generalizedzykovstheorem,
      title={{Generalized Zykov's Theorem}}, 
      author={Rajat Adak and L. Sunil Chandran},
    journal={arXiv preprint arXiv:2512.02958},
      year={2025},
      doi={10.48550/arXiv.2512.02958}
}

@article{adak2025vertexbasedlocalizationgeneralizedturan,
      title={{Vertex-Based Localization of generalized Tur\'{a}n Problems}},
      author={Adak, Rajat and Chandran, L Sunil},
      journal={arXiv preprint arXiv:2508.20936},
      year={2025},
      doi={10.48550/arXiv.2508.20936}
}

\end{document}